\documentclass[a4paper,12pt]{report}
\usepackage{amsmath,amssymb}
\usepackage[brazil]{babel}
\usepackage[latin1]{inputenc}
\usepackage{graphicx}
\numberwithin{equation}{section}

\newtheorem{teo}{Theorem}

\newtheorem{cory}{Corollary}

\newcommand{\whitebox}{$\Box$}

\newenvironment{proof*}[1][Proof]{\textbf{#1:} }

\newcommand{\NN}{\mathbb{N}}
\newcommand{\ZZ}{\mathbb{Z}}

\def\CaixaPreta{\vrule Depth0pt height5pt width5pt}
\def\comecaprova{\noindent {\bf Proof:} \hspace{2mm}}
\def\terminaprova{\hfill \CaixaPreta \vspace{3mm}}

\newcommand{\mapdown}[2]{{#1}\kern -0.05cm\raise 0.02cm
{#2}}

\begin{document}

\begin{center}
{\LARGE On the definability of rational integers in a class of rings}\vspace{0,5cm}\\ {\large Eudes Naziazeno}\\Universidade Federal de Pernambuco\vspace{0.5cm}
\end{center}
{\bf Abstract:} We provide a sufficient condition for a ring to have a first-order definition for $\mathbb{Z}$.\vspace{0.2cm}

\section*{Section 0: Introduction.} Raphael Robinson in [1] presented a list of undecidable rings. In Section 2, Robinson showed that if $R$ is a commutative integral domain of characteristic zero and $R$ is definable in $R[x]$, then $\mathbb{Z}$ is definable in $R[x]$. He used field extensions to prove this result.\vspace{0.2cm}\\In this article we provide a sufficient condition (more elementary than Robinson's) for a ring $R$ to have $\mathbb{Z}$ as a definable set. As a consequence, Robinson's result can be viewed as a corollary.

\section*{Section 1: The main theorem and its proof} The general theorem of this work is the following:
\begin{teo}\label{teo1}
Let $R$ be a ring (not necessarily commutative) such that $\ZZ\subseteq R$. Suppose that there is a definable set $A$ such that $\ZZ\subseteq A\subseteq R$ and that for all $a\in A$ and $n>0$ integer we have that $n-a\in A$. Suppose also that $R$ has an element $p$ such that:
\begin{enumerate}
  \item $\operatorname{POW}(p):=\{p,p^2,p^3,p^4,\cdots\}$ is an infinite definable set (possibly using $p$ as a parameter) and all $p^n$ are pairwise different;
  \item $p-1$ is not a zero divisor; and
  \item $\forall a\in A\ \ (p-1)|a \Rightarrow a=0$.
\end{enumerate}
Then $\ZZ$ is definable in $R$.
\end{teo}
\comecaprova Let $\psi(y)$ be a definition for $\operatorname{POW}(p)$, i.e., $\psi(y)\Leftrightarrow y=p^n$ for some $n\in\NN$. Consider the sentence
$$\varphi(t)\ :\ \exists y\ \ \psi(y)\wedge\ p^2|y\wedge\ \exists w\ \ y-1=w\cdot(p-1)\wedge\ t\in A\wedge (p-1)|(w-t).$$
We claim that $z\in\NN\Leftrightarrow z=0\vee z=1\vee\varphi(z)$. Indeed, suppose that $z\neq0$ but $\varphi(z)$. Then, $y=p^n$ for some $n\geq2$ integer. Observe that
$$y-1=p^n-1=(p^{n-1}+\cdots+p+1)\cdot(p-1)$$
Since $\exists w\ \ y-1=w\cdot(p-1)$, then
$$w\cdot(p-1)=(p^{n-1}+\cdots+p+1)\cdot(p-1)$$
$$(w-(p^{n-1}+\cdots+p+1))\cdot(p-1)=0$$
Since $p-1$ is not a zero divisor, we have that $w=p^{n-1}+\cdots+p+1$. For all $t\in R$ we can write
$$w-t=(p^{n-2}+2p^{n-3}+\cdots+(n-2)p+(n-1))\cdot(p-1)+n-t$$
So, $(p-1)|(w-t)\Rightarrow(p-1)|(n-t)$. Since $n\in\ZZ$ and $t\in A$, then $n-t\in A$. Therefore, $(p-1)|(n-t)\Rightarrow t\in\ZZ$. The converse is easy, because for each $z=n\geq2$ we can set $y=p^n$.\terminaprova\vspace{0.2cm}\\
Concluding this section, we emphasize that we reduced the problem of defining $\ZZ$ in $R$ to the problem of finding that special set $A\subseteq R$ and that special element $p\in R$. On the next section we provide two cases for which it is not hard to find these special tools.

\section*{Section 2: Commutative and noncommutative examples.}
The first corollary is the same as in $[1]$, and the other two are new ones.
\begin{cory}\label{cory1}
Let $R[x]$ be a commutative integral domain such that $\operatorname{char} R=0$ and $R$ is definable. Then, $\ZZ$ in definable in $R[x]$.
\end{cory}
\comecaprova
We claim that $A:=R$ and $p:=x$ satisfy the conditions on Theorem \ref{teo1}. Indeed, trivially $x-1$ is not a zero divisor. Suppose that $f(x)=a_0+a_1x+\cdots+a_nx^n$ is such that $a=(x-1)\cdot f(x)$ and $a\in R$. Then,
$$a=(x-1)\cdot(a_0+a_1x+\cdots+a_nx^n)$$
$$a=a_0x+a_1x^2+\cdots+a_nx^{n+1}-a_0-a_1x-\cdots-a_nx^n$$
$$a=-a_0+(a_0-a_1)x+\cdots+(a_{n-1}-a_n)x+a_nx^n$$
Since these polynomials are equal, then $a=a_0=a_1=\cdots=a_n=0$. Then, $\forall a\in R\ \ (x-1)|a \Rightarrow a=0$. Now, consider the sentence
$$\phi(t):\ \forall d\ \ (d\nmid1\wedge d|t)\Rightarrow x|d\wedge\ (x-1)|(t-1)$$
We claim that $z\in\operatorname{POW}(x)\Leftrightarrow\phi(z)$. Indeed, let us first note that $x$ is prime in $R[x]$ whenever $R$ is an integral domain.  This implies that
$$(\forall d\ \ (d\nmid1\wedge d|t)\Rightarrow x|d)\Leftrightarrow z=ax^n$$
for which $a$ is invertible. Since $R$ is an integral domain, then $a\in R$. Now we prove that $\phi(z)\Rightarrow z\in\operatorname{POW}(p)$ (the converse is trivial). Suppose that $\phi(z)$. If $z$ is invertible, we are done (just take $z=ax^0$). Now suppose that $z$ is not invertible. Since $z|z$ and $z\nmid1$, then $x|z$. So, $z=z_1\cdot x$. If $z_1|1$, we are done. If $z_1\nmid1$, then $z_1=z_2\cdot x$, for some $z_2$. Then, $z=z_2\cdot x^2$. We can go on with this process until we find some $a$ such that $z=a\cdot x^n$ and $a|1$. This process has to stop because $z\in R[x]$ has finite degree and $R$ is an integral domain. Observe that
$$z-1=ax^n-1=a(x^n-1)+a-1=a(x^{n-1}+\cdots+x+1)\cdot(x-1)+a-1$$
So, $(x-1)|(z-1)\Rightarrow (x-1)|(a-1)$. Since $a|1$ and $R$ is an integral domain, then $a-1\in R$. By our last claim, $a=1$, and so $z=x^n$. Therefore, with this we proved that $A:=R$ and $p:=x$ satisfy all the hypothesis on Theorem \ref{teo1}. Thus, $\ZZ$ is definable in $R[x]$.\terminaprova\vspace{0.2cm}\\Robinson used field extensions to prove this last result. However, as the reader can see, our proof consists of arguments that are more elementary than Robinson's. As another corollary, we see (as Robinson pointed in [1]), that if $R$ is a commutative field with $\operatorname{char} R=0$, then $\ZZ$ is definable in $R[x]$\ (because $z\in R\Leftrightarrow z=0\vee z|1$). Now we present a similar result, but over a noncommutative ring.\vspace{0.2cm}\\The {\it quantum affine plane over R}, denoted by $R_q[x,y]$, is the noncommutative ring of the polynomials with complex coefficients, on the variables $x$ and $y$, such that $$xy=qyx$$
for some fixed $q\in R^{*}$ (see [2]). For example, in $\mathbb{C}_2[x,y]$ we have that $(3+x)\cdot(2+y)=6+2x+3y+xy$, but $(2+y)\cdot(3+x)=6+2x+3y+2xy\ .$
\begin{cory}
$\ZZ$ is definable in $\mathbb{C}_q[x,y]$.
\end{cory}
\comecaprova
With $A:=\mathbb{C}$ and $p:=x$, we can show that $A$ and $p$ satisfy all the conditions of Theorem \ref{teo1}. The proof of this fact is basically the same as in Corollary \ref{cory1}. However, the only thing we have to adjust is that if we write
$$f|g\Leftrightarrow\exists h\ \ g=f\cdot h$$
then we have to write $x^n-1=(x-1)\cdot(x^{n-1}+\cdots+x+1)$, i.e., with $x-1$ on the left side instead of of writing as we did before (and noticing that $z\in\mathbb{C}\Leftrightarrow z=0\vee z|1$). All the calculations are similar to those on Corollary \ref{cory1} (the only difference is that we have to deal with powers of $q$ over the equations), and that is why we avoided to repeat it here.\terminaprova\vspace{0.2cm}\\ Finally, we end with a conjecture and one of its corollaries.\vspace{0.5cm}\\{\bf Conjecture:} {\it If $R$ is a reduced commutative ring, then $\operatorname{POW}(x)$ is definable in $R[x]$ and in $R_q[x,y]$, possibly making use of a formula for $\ulcorner z\in R\urcorner$.}\vspace{0.2cm}

\begin{cory}
Let $R$ be a reduced commutative ring with $\operatorname{char}R=0$. If $R$ is definable in $R[x]$, then $\ZZ$ is definable in $R[x]$. If $R$ is definable in $R_q[x,y]$, then $\ZZ$ is definable in $R_q[x,y]$.
\end{cory}
\comecaprova
We can use adaptations of the corollaries before, noticing that since $R$ has no nilpotent elements, then the invertibles in $R[x]$ and in $R_q[x,y]$ are precisely the invertibles of $R$ (see [3]).\terminaprova

\section*{Acknowledgements}
We would like to express our sincere thanks and appreciation to Thomas Scanlon, Alexandra Shlapentokh, and Carlos Videla, for their kindness and the inspiring advices.

\section*{References}

[1] Robinson, R. M. {\it Undecidable Rings}, Trans. Amer. Math. Soc. 70 (1951), 137-159.\vspace{0.5cm}\newline
[2]K. R. GOODEARL and E. S. LETZTER {\it Quantum $n$-space as a quotient of classical $n$-space} Trans. Amer. Math. Soc. 352 (2000), 5855-5876.\vspace{0.5cm}\newline
[3] M. F. ATIYAH and I. G. MacDONALD {\it Introduction to Commutative Algebra}, Addison-Wesley (1969).

\end{document}